\documentclass[a4paper,reqno]{amsart}
\usepackage{amssymb}
\usepackage{latexsym}
\usepackage{amsmath}
\usepackage{euscript}
\usepackage{graphics,color}
\usepackage[all]{xy}
\usepackage[margin=3cm]{geometry}
\usepackage{MnSymbol} 

\newcommand{\be}{\begin{equation}}
\newcommand{\ee}{\end{equation}}
\newcommand{\ba}{\begin{eqnarray}}
\newcommand{\ea}{\end{eqnarray}}
\newcommand{\baa}{\begin{eqnarray*}}
\newcommand{\eaa}{\end{eqnarray*}}
\newcommand{\bb}{\begin{thebibliography}}
\newcommand{\eb}{\end{thebibliography}}

\newcommand{\bi}[1]{\bibitem{#1}}
\newcommand{\lab}[1]{\label{#1}}
\newcommand{\re}[1]{(\ref{#1})}

\newcommand*{\fdots}{\genfrac{}{}{0pt}{}{}{\cdots}}
\newcommand*{\fminus}{\genfrac{}{}{0pt}{}{}{-}}

\newcounter{my}
\newcommand{\he}%
   {\stepcounter{equation}\setcounter{my}%
   {\value{equation}}\setcounter{equation}0%
   }%
\newcommand{\she}%
   {\setcounter{equation}{\value{my}}%
    }%

\renewcommand\t{\tilde}

\newcommand\vphi{\varphi}

\newtheorem{pr}{Proposition}

\theoremstyle{definition}

\numberwithin{equation}{section}

\begin{document}

\title{Classical Sturmian sequences}

\author{Alexei Zhedanov}

\address{School of Mathematics, Renmin University of China, Beijing 100872,CHINA}

\begin{abstract}
The Sturm sequence is generated by a pair of polynomials $P(x)$ and $P'(x)$, where $P(x)$ is assumed to have simple real roots. Euclidean algorithm  generates then a finite sequence of polynomials orthogonal on the grid $x_s$ of roots of the polynomial $P(x)$. This algorithm can be exploited in order to find the number of roots of the polynomial $P(x)$ inside a given interval. We consider the "inverse" problem: what is the explicit system of orthogonal polynomials generated by the prescribed grid $x_s$ of "classical" type. The main results are the following. The generic linear grid generates a special case of the Hahn polynomials. The quadratic grids $x_s=x(s+1)$ and $x_s=s(s+2)$ correspond to two special cases of the Racah polynomials. The generic exponential grid is related to a special case of the q-Hahn polynomials. Finally, we show that two special trigonometric grids are related to the Chebyshev polynomials of the first and second kind.   
\end{abstract}

\keywords{}


\maketitle

\section{Introduction}
\setcounter{equation}{0}
The classical Sturm algorithm \cite{Waerden} allows to find the number of real roots of the polynomial $P(x)$ inside an interval. The main idea of the algorithm is application the Euclidean algorithm which generates a sequence of polynomials of decreasing degrees starting from the pair of polynomials $P(x)$ and $P'(x)$ where $P'(x)$ is the derivative of $P(x)$. Assume for simplicity that all roots of the polynomial $P(x)$ are real and simple. Then all the roots of the polynomial $P'(x)$ are also real and simple and interlace the roots of the polynomial $P(x)$. This means that every root of the polynomial $P'(x)$ is located between two neighbor roots of the polynomial $P(x)$.

In more general situation, assume that $P_{N+1}(x)$ is a monic polynomial with $N+1$ simple real zeros $x_s$. These zeros are ordered by increasing
\be
x_0 <x_1 <x_2 < \dots x_N . \lab{order_x} \ee
Assume also the monic polynomial $P_N(x)$ of degree $N$ is arbitrary with the only condition that its simple real zeros interlace zeros of $P_{N+1}(x)$. Dividing $P_{N+1}(x)$ by $P_N(x)$ we get
\be
P_{N+1}(x) = (x-b_N)P_N(x) - u_N P_{N-1}(x), \lab{rec_N+1} \ee
where $P_{N-1}(x)$ is a monic polynomial of degree $N-1$ with zeros interlacing zeros of $P_{N}(x)$. The coefficient $u_N$ is strictly positive $u_N>0$. This process can be continued  which yields the chain of monic polynomials $P_{N-2}(x), P_{N-3}(x), \dots, P_0=1$. They satisfy the recurrence relation
\be
P_{n+1}(x) = (x-b_n)P_n(x) - u_n P_{n-1}(x), \quad n=1,2, \dots, N \lab{rec_n} \ee 
with strictly positive coefficients $u_n>0$. 

On the other hand, recurrence relation \re{rec_n} (together with the condition $u_n>0$) is equivalent to the statement that polynomials $P_n(x), \: n=0,1,\dots, N$ are orthogonal
\be
\sum_{s=0}^N w_s P_n(x_s) P_m(x_s) = h_n \: \delta_{nm} \lab{ort_P} \ee
where the normalization coefficient 
\be
h_n = u_1 u_2 \dots u_n, \quad h_0 =1, \lab{h_n} \ee
and where the positive weights $w_s>0$ are  \cite{Chi} 
\be
w_s= \frac{h_N}{P_{N+1}'(x_s) P_N(x_s)} . \lab{w_s_gen} \ee
The weights $w_s$ are normalized
\be
\sum_{s=0}^N w_s =1. \lab{norm_w} \ee
Thus starting with any pair $P_{N+1}(x), P_N(x)$ of monic polynomials with simple interlacing zeros, one can reconstruct a family of polynomials $P_n(x), n=0,1, \dots, N-1$ which are orthogonal on zeros $x_s$ of the polynomial $P_{N+1}(x)$ with positive weight function. Hence, the Jacobi matrix 
\[
J =
 \begin{pmatrix}
  b_{0} & 1 & 0 &    \\
  u_{1} & b_{1} & 1 & 0  \\
   0  &  u_2 & b_2 & 1    \\
   &   &  \ddots &    \ddots  \\
      & & \dots &u_{N-1} & b_{N-1} & 1 \\   
            & & \dots &0 & u_N & b_N
         
\end{pmatrix}.
\] 
can uniquely be associated with the pair $P_{N+1}(x), P_N(x)$ \cite{Chi}.

The mirror-dual Jacobi matrix $J^*$ has the entries
\be
b_n^*=b_{N-n}, \quad u_n^*=u_{N+1-n} , \quad n=0,1,\dots N \lab{ub*} \ee
(it is assumed in \re{ub*} that $u_0=u_{N+1}=0$). Algebraically, the mirror dual matrix $J^*$ satisfies the relation
\be
J^* = RJ^TR, \lab{RJ} \ee
where $J^T$ is transposed Jacobi matrix and where $R$ is the reflection matrix  \cite{persym}
\[
R=\begin{pmatrix}
  0 & 0 & \dots & 0 & 1    \\
  0 & 0 & \dots  & 1 & 0  \\
    \dots  & \dots & \dots & \dots & \dots      \\
   1 & 0 &  \dots & 0 &0  \\
\end{pmatrix}.
\] 
Corresponding mirror-dual orthogonal polynomials $P_n^*(x)$ are uniquely defined by the recurrence relation
\be
P_{n+1}^*(x) = (x-b_n^*)P_n^*(x) - u_n^* P_{n-1}^*(x) \lab{rec_dual} \ee
and by the initial conditions
\be
P_0^*=1, \; P_1^*(x) = x-b_0^*=x-b_N . \lab{ini_dual} \ee
The mirror-dual polynomials satisfy the orthogonality relation
\be
\sum_{s=0}^N w_s^* P_n^*(x_s) P_m^*(x_s) = h_n^* \: \delta_{nm}, \lab{dual_ort_P} \ee
where the (normalized) the dual weights $w_s^*$ are
\be
w_s^* = \frac{P_N(x_s)}{P_{N+1}'(x_s)}. \lab{w*} \ee
Comparing \re{w*} with \re{w_s_gen} we arrive at the important duality relation between the weights \cite{BS}, \cite{Borodin}, \cite{VZ_dual}
\be
w_s w_s^* = \frac{h_N}{(P_{N+1}'(x_s))^2}. \lab{ww_rel} \ee
So far, the polynomial $P_N(x)$ was an arbitrary monic polynomial of degree $N$ with the only condition of interlacing zeros with respect to $P_{N+1}(x)$. 
The Jacobi matrices $J$, $J^*$ and orthogonal polynomials $P_n(x)$, $P_n^*(x)$ are in one-to-one correspondence with the choice of the polynomials $P_N(x), P_{N+1}(x)$. 
Consider now the special Sturmian case when $P_N(x) = (N+1)^{-1} P_{N+1}'(x)$. Then from \re{w*} we get
\be
w_s^* = \frac{1}{N+1} . \lab{w_s_const} \ee
Thus for the Sturmian case all the dual weights $w_s^*$ are the same (i.e. they do not depend on $s$). But the weights $w_s^*$ and the grid $x_s$ define uniquely the matrices $J$ and $J^*$ and hence - the pair of polynomials $P_{N+1}(x), P_{N}(x)$. It is natural to call the polynomials with constant weights \re{w_s_const} the {\it Legendre-type} polynomials. Indeed, the Legendre polynomials are orthogonal with the constant weight $1/2$ on the interval $[-1,1]$ \cite{KLS}. 

We thus have 
\begin{pr}
The Sturmian pair $P_{N+1}(s)$ and $P_N(x) = (N+1)^{-1} P_{N+1}'(x)$ corresponds to mirror-dual orthogonal polynomials $P_n^*(x)$ of Legendre type. Conversely, any finite system of orthogonal polynomials $P_n^*(x), \: n=0,1,\dots, N$ with the weights $w_s^*=(N+1)^{-1}$, corresponds to the Sturmian pair $P_{N+1}(x)$ and $P_N(x) = (N+1)^{-1} P_{N+1}'(x)$.
\end{pr}  
We can therefore reduce the problem of classification of orthogonal polynomials $P_n(x)$ corresponding to the Sturmian pair $P_{N+1}(s)$ and $P_N(x) = (N+1)^{-1} P_{N+1}'(x)$, to the problem of finding polynomials $P_n^*(x)$ of Legendre type orthogonal with respect to constant weights $w_s^*=(N+1)^{-1}$ on the grid $x_s$ of zeros of the polynomial $P_{N+1}(x)$.

Note that earlier similar property was established with respect to classical orthogonal polynomials \cite{DS}, \cite{VZ_dual}. Indeed, consider, e.g. the Hermite polynomials  $H_n(x)$ \cite{KLS}. They satisfy the relation
\be
H_{n+1}'(x) = (n+1) H_n(x), \quad n=0,1,2,\dots .  \lab{Her_diff} \ee 
Hence for any $N$ the pair $H_{N+1}(x), H_N(x)$ is the Sturmian pair. In turn, it follows that for any $N$ the mirror-dual polynomials $H_n^*(x)$ are orthogonal with constant weights $(N+1)^{-1}$ on the grid $x_0, x_1, \dots, x_N$ of the roots of the Hermite polynomial $H_{N+1}(x)$ \cite{DS}, \cite{VZ_dual}.

In this case however, the grid $x_s$ of the roots of the Hermite polynomial $H_{N+1}(x)$ is not "classical": there is no explicit expression for these roots as a function of $s$. Instead, we consider the following problem. Let us start with prescribed "classical" grid $x_s$ with distinct real entries. This grid defines uniquely the monic polynomial $P_{N+1}(x)$. Consider then the monic polynomial $P_N(x) = (N+1)^{-1} P_{N+1}'(x)$ of degree $N$. Clearly, the roots of polynomials $P_N(x)$ and $P_{N+1}(x)$ are interlacing.  Together, the polynomials $P_{N+1}(x), P_N(x)$ constitute the Sturmian pair.   Our main problem will be deriving explicit expressions for the Jacobi matrices $J$ (or $J^*$) and for corresponding orthogonal polynomials $P_n(x)$. 

By "classical grid" we mean the orthogonality grids for classical orthogonal polynomials of discrete variable \cite{KLS}. The most general grid of such type should satisfy the linear difference equation   \cite{NSU},  \cite{Ter}, \cite{VZ_Bochner}
\be
x_{s+1} + x_{s-1} - \Omega x_s = \nu \lab{eq_grid} \ee
with some real constants $\Omega$ and $\nu$. The type of the grid depends on the parameter $\Omega$. More exactly:

\vspace{2mm}

(ia) if $|\Omega| >2$ then one has the Askey-Wilson grid
\be
x_s = C_1 q^s + C_2 q^{-s} + C_0, \lab{AW_grid} \ee
where $\Omega = q+q^{-1}$ with some real $q$. In particular, when $C_1=0$ we have the exponential grid 
\be
x_s=C_2 q^{-s} + C_0 . \lab{exp_grid} \ee

\vspace{2mm}

(ib) the case $|\Omega| <2$ corresponds to trigonometric Askey-Wilson grid \re{AW_grid} with $q=\exp(i \theta)$ (see \cite{SZ_Duke} for details).

\vspace{2mm}

(ii) if $\Omega =2$ then one has the quadratic grid
\be
x_s = C_2 s^2 + C_1 s + C_0 . \lab{quad_grid} \ee
In particular, when $C_2=0$ we have the linear grid 
\be
x_s = C_1 s + C_0. \lab{l_grid} \ee

\vspace{2mm}

(iii) if $\Omega=-2$ then one has the Bannai-Ito grid
\be
x_s = (-1)^s \left( C_1 s + C_2\right) + C_0. \lab{BI_grid} \ee

We restrict ourselves only with linear, exponential and quadratic grids and with a special type of trigonometric grid. Other types of grids will be considered elsewhere.

The paper is organized as follows. In Section 2, we consider generic properties of polynomials $P_n(x), P_n^*(x)$ corresponding to the Sturm sequence. In particular, we reinterpret Sylvester's approach to the Sturm algorithm in terms of orthogonal polynomials. In section 3, the linear grid $x_s=s$ is considered. In this case the mirror duals with respect to the Chebyshev-Legendre polynomials appear as general solution. In Section 4, we study the quadratic grid $x_s=s(s+\tau)$, where $\tau$ is the only essential parameter which characterizes the grid. For the two cases: $\tau=1$ and $\tau=2$ we obtain solutions in terms of Racah polynomials.  In Section 5, the exponential grid $x_s=q^{-s}$ is analyzed leading to the general solution related to the q-Hahn polynomials. In Section 6, two types of trigonometric grids are considered. They correspond  to solutions in terms of Chebyshev polynomials of the first and second kind. In concluding remarks in Section 7 we discuss some open problems.

\section{Generic properties of Sturmian sequences}
\setcounter{equation}{0}
Before considering concrete examples, let us mention some generic properties of the Jacobi matrices $J, J^*$ and corresponding polynomials $P_n(x), P_n^*(x)$. First of all, the affine transformation  of the grid $x_s \to \alpha x_s + \beta$ leads to the polynomials $P_n(x)$ with corresponding affine transformed argument $x$. For the Jacobi matrices $J, J^*$ this transformation leads to shifting of the diagonal entries $b_n \to b_n + const$ and to rescaling the off-diagonal entries $u_n$.    Hence we can use this affine transformation in order to reduce the grid to the most convenient form. 

Second, let us consider the moments $c_n^*$ corresponding to the dual weights $w_s^*$
\be 
c_n^* = \sum_{s=0}^N w_s^* x_s^n = (N+1)^{-1} \sum_{s=0}^N  x_s^n . \lab{moms*} \ee
This formula indicates that the moments $c_n^*$ are nothing else than sum of powers of the grid points. For the classical grids these sums can be explicitly calculated. For example, for the linear grid (after its appropriate affine transformation) $x_s=s, \: s=0,1,2,\dots, N$ we have that the moments $c_n^*$ coincide with Faullhaber sums 
\be
c_n^*=\sum_{s=0}^N s^n .\lab{c_F} \ee
These sums can be expressed in terms of the Bernoulli numbers  (see e.g. \cite{Beardon}).

We can relate the above results with Sylvester's approach to the Sturm problem \cite{Syl}, \cite{Prasolov}.

Indeed, let $P_n(x)$ be a finite system of orthogonal polynomials defined by the recurrence relation \re{rec_n}. One can introduce then the system of orthogonal polynomials of second kind $R_n(x)$  defined by \cite{Chi}
\be 
R_{n+1}(x) + b_{n+1}R_n(x) + u_{n+1}(x) = xR_n(x), \quad n=1,2,\dots, N, \quad R_0 =1, \; R_1(x) = x-b_1 .\lab{R_rec} \ee
Then one can construct the continued fraction \cite{Chi}
\begin{equation}
F(z) = \frac{R_{N}(z)}{P_{N+1}(z) } =
\frac{1}{z-b_0}\fminus \frac{u_1}{z-b_1}\fminus
\frac{u_2}{z-b_2}
\fminus\fdots \fminus \frac{u_N}{z-b_N}
 .\label{F_CF}
\end{equation}
The rational function $F(z)$ plays the role of the Stieltjes transform of the orthogonality measure. Indeed, expanding $F(z)$ in term of partial fractions, we have \cite{Chi}
\be
F(z) = \sum_{s=0}^N \frac{w_s}{z-x_s} \lab{F_par} \ee
so that the weights $w_s$ are residues of $F(z)$ at simple poles $z=x_s$.

On the other hand, by the Euclidean algorithm, we have the continued fractions
\be
F^*(z)=\frac{P_{N}(z)}{P_{N+1}(z) } =
\frac{1}{z-b_N}\fminus \frac{u_N}{z-b_{N-1}}\fminus
\frac{u_{N-1}}{z-b_{N-2}}
\fminus\fdots \fminus \frac{u_1}{z-b_0}
. \label{F_*} \ee
It is clear that the functions $F(z)$ and $F^*(z)$ are mirror conjugate. In more details, the function $F(z)$ corresponds to the Jacobi matrix $J$ while the function $F^*(z)$ corresponds to the Jacobi matrix $J^*$. By definition, the simple poles $x_s$ are the same for both functions $F(z)$ and $F^*(z)$. Hence, if we need only information about poles $x_s$, we can use any of the Stieltjes functions $F(z)$ or $F^*(z)$. But for the Sturmian pair $P_{N}(x)= (N+1)^{-1}P_{N+1}'(x)$, the function $F^*(z)$ has extremely simple expansion in terms of partial fraction:
\be
F^*(z) = \frac{P_{N+1}'(z)}{(N+1)P_{N+1}(z)} = \frac{1}{N+1} \: \sum_{s=0}^N \frac{1}{(z-x_s)} . \lab{F*_par} \ee
In other words, formula \re{F*_par} is equivalent to the statement that the "mirror" weights $w_s^*$ are all equal one to another \re{w_s_const}.  And this was the key idea of Sylvester's approach - to use the mirror-reflected continued fraction \re{F_*} instead of \re{F_CF} \cite{Syl}. We thus see that Sylvester's approach has simple and transparent interpretation in terms of orthogonal polynomials.

As a direct application, let us consider the Hankel determinants constructed from abstract moments $c_n$ 
\be
\Delta_n = det \begin{pmatrix}
  c_0 & c_1 & \dots & c_{n-1} & c_n   \\
  c_1 & c_2 & \dots  & c_n & c_{n+1}  \\
    \dots  & \dots & \dots & \dots & \dots      \\
   c_n & c_{n+1} &  \dots & c_{2n-1} & c_{2n}  \\
\end{pmatrix}, \quad \Delta_n^{(1)}= \begin{pmatrix}
  c_1 & c_2 & \dots & c_{n} & c_{n+1}   \\
  c_2 & c_3 & \dots  & c_{n+1} & c_{n+2}  \\
    \dots  & \dots & \dots & \dots & \dots      \\
   c_{n+1} & c_{n+2} &  \dots & c_{2n} & c_{2n+1} \\
\end{pmatrix}.\lab{Hankels} \ee 
The classical result from the theory of orthogonal polynomials is the Stieltjes criterion \cite{Chi}: the conditions
\be
\Delta_n >0, \quad \Delta_n^{(1)} >0 \quad n=0,1,2,\dots \lab{Stielt_crit} \ee
are equivalent to existence of a positive orthogonality measure on the semi-axis $0<x<\infty$. This means that there exists a positive measure $d \mu(x)$ such that
\be
c_n = \int_{0}^{\infty} x^n d \mu(x), \quad n=0,1,2,\dots \lab{c_n_def} \ee
In our case this criterion means that all roots of the polynomial $P_{N+1}(x)$ are positive. The above criterion remains true if one replaces the  moments $c_n$ with the "mirror" moments $c_n^*$ defined by \re{moms*}. This is one of the main Sylvester's  results \cite{Prasolov}.

Finally, note that if the grid $x_s$ is symmetric around zero, i.e. if
\be
x_{N-s}=-x_s, \quad s=0,1, \dots N, \lab{sym_x} \ee
then the matrices $J,J^*$ have zero diagonal entries, i.e. $b_n=0, \: n=0,1,\dots,N$. Equivalently, this means that the polynomials $P_n(x)$ and $P_n^*(x)$ are symmetric:
\be
P_n(-x) = (-1)^n P_n(x) . \lab{sym_P} \ee 
This property can easily be derived from the observation that for symmetric grid $x_s$  \re{sym_x} the characteristic polynomial $P_{N+1}(x)$ is either odd or even depending on parity of $N$. The polynomial $P_N(x)$ has obviously opposite parity. Then by induction one can prove that all further polynomials $P_n(x)$ satisfy the property \re{sym_P}. Equivalently, this means that the Jacobi matrices $J$ and $J^*$  have zero diagonal entries \cite{Chi}.

\section{The linear grid}
\setcounter{equation}{0}
We start with the most simple case when the grid is linear. By an affine transformation it is always possible to reduce the linear grid to the standard form
\be
x_s =0,1,2,\dots, N . \lab{lin_grid} \ee
Then we arrive at the problem of finding orthogonal polynomials having constant weights $w_s^* =(N+1)^{-1}$ on the linear grid \re{lin_grid}. This problem was solved by Chebyshev \cite{NSU}. It appears that these Chebyshev polynomials of Legendre type are special case of the Hahn polynomials (Note nevertheless, that the generic Hahn polynomials were introduced by Chebushev too \cite{KLS}).

In order to describe these Chebyshev-Legendre polynomials, recall first that generic monic Hahn polynomial $H_n(x;\alpha,\beta,N)$ depend on 2 parameters $\alpha, \beta$ and on the integer parameter $N$. Explicitly, these polynomials are expressed in terms of hypergeometric function
\be
H_n(x; \alpha, \beta,N) = \kappa_n \:{_3}F_2 \left( {-n, -x, n+\alpha+\beta+1 \atop -N, \alpha+1}; 1 \right) ,  \lab{Hahn_H} \ee
where $\kappa_n$ is the normalization coefficient 
\[
\kappa_n = \frac{(-N)_n (\alpha+1)_n}{(n+\alpha+\beta+1)_n}
\]
and where $(x)_n=x(x+1) \dots (x+n-1)$ is the Pochhammer symbol.

The Hahn polynomials satisfy the three-term recurrence relation 
\be
H_{n+1}(x) + b_n H_n(x) + u_n H_{n-1}(x) = x H_n(x) \lab{3-Hahn} \ee
with 
\be
b_n = A_n+C_n, \quad u_n = A_{n-1} C_n, \lab{bu_H} \ee
where
\be
A_n = \frac{(n+\alpha+\beta+1)(n+\alpha+1)(N-n)}{(2n+\alpha+\beta+1)(2n+\alpha+\beta+2)}, \quad  C_n = \frac{n(n+\alpha+\beta+N+1)(n+\beta)}{(2n+\alpha+\beta+1)(2n+\alpha+\beta)}. \lab{AC_H} \ee
These polynomials are orthogonal on the linear lattice $x_s=s$
\be
\sum_{s=0}^N W_s H_n(s;\alpha,\beta,N) H_m(s;\alpha,\beta,N) = h_n \delta_{nm} \lab{ort_H} \ee
with the weights
\be
W_s = \nu \: \frac{(-N)_s (\alpha+1)_s }{s!(-N-\beta)_s}, \lab{W_s} \ee
where 
\be
\nu = \frac{(\beta+1)_N}{(\alpha+\beta+2)_N} \lab{nu_W} \ee
is the normalization coefficient needed to fulfill the standard condition
\be
\sum_{s=0}^N W_s =1 . \lab{norm_W} \ee 
It is now easy to see that if $\alpha=\beta=0$ then 
\be
W_s =\frac{N!}{(N+1)!} = \frac{1}{N+1} \lab{W_s_Cheb} \ee
and we arrive at the Chebyshev polynomials of Legendre type with constant wrights on the set of points $x_s=s$.

Hence we have for the mirror-dual polynomials $P_n^*(x)$ the recurrence coefficients 
\be
u_n^*= \frac{n^2 ((N+1)^2-n^2)}{4(4n^2-1)} , \quad b_n^* =N/2, \quad n=0,1, \dots, N \lab{ub_mirr_uni} \ee
The corresponding Sturmian polynomials $P_n(x)$ also belong to the  class of Hahn polynomials. Indeed, it is easy to check that the transformation of the parameters
\be
\t \alpha = -N-1 - \beta, \quad \t \beta = -N-1-\alpha \lab{mirror_ab_H} \ee
is equivalent to the mirror transform $u_n \to u_n^*=u_{N+1-n}, \: b_n \to b_n^*=b_{N-n}$ of the recurrence coefficients of the Hahn polynomials. Hence
\be
H_n^*(x; \alpha, \beta,N) = H_n(x, \t \alpha, \t \beta,N) \lab{tilde_H} \ee    
For the case of Chebyshev polynomials of Legendre type we have
\be
\tilde \alpha = \tilde \beta = -N-1 \lab{ab_Sturm_H} \ee

We thus arrive at the 
\begin{pr} 
The Sturmian sequence of the polynomials $P_n(x)$ on the linear grid $x_s=s=0,1,\dots,N$ coincides with the set of the Hahn polynomials 
\be
P_n(x) = H_n(x, \alpha, \beta, N), \lab{P_H_Sturm} \ee
where $\alpha=\beta=-N-1$.
\end{pr}
It is interesting to note that polynomials \re{P_H_Sturm} are orthogonal on the grid $x_s=s$ with respect to the squared binomial distribution
\be
w_s = \nu \: \left( \frac{N!}{s!(N-s)!}\right)^2, \lab{sq_binom} \ee
where
\be
\nu = \frac{(N!)^2}{(2N)!} \lab{nu_St_H} \ee
Moreover, these polynomials satisfy the simple difference equation on the linear grid
\be
(x-N)^2 \left(P_n(x+1)-P_n(x) \right) + x^2 \left(P_n(x-1)-P_n(x) \right) = n(n-2N-1) P_n(x) . \lab{dfr_Hahn } \ee

\section{Quadratic grid}
\setcounter{equation}{0}
Consider the quadratic grid $x_s = a_2s^2 + a_1 s + a_0$, where all parameters $a_0,a_1,a_2$ are real.  By an appropriate affine transformation it is always possible to reduce this grid to the form
\be
x_s = s(s+\tau) \lab{quad_grid_red} \ee
with the only parameter $\tau$ which defines the type of the grid. Note that for the linear grid there are no free parameters: all linear grids are affine equivalent to the standard one $x_s=s$. For the quadratic grid we should distinguish different types of grids depending on the essential parameter $\tau$.

As in the previous section, we should first to construct the mirror-dual polynomials $P_m^*(x)$ of Legendre type on the grid \re{quad_grid_red}. This means that these polynomials should  satisfy the  orthogonality relation
\be
\sum_{s=0}^N  P_n^*(s(s+\tau)) P_m^*(s(s+\tau)) = (N+1) h_n^* \: \delta_{nm} \lab{quad_mirror_ort} \ee
It is natural to search the polynomials $P_n^*(x)$ among "classical" polynomials orthogonal on the quadratic grid \re{quad_grid_red}. Such polynomials are well known. The most general family of these polynomials - Racah polynomials - contains 4 parameters $\alpha,\beta, \gamma, \delta$ \cite{KLS}. One of these parameters, say $\alpha$, should be chosen as
\be
\alpha= -N-1 \lab{alpha_N+1} \ee  
to fulfill the truncation condition. The explicit expression of these polynomials is \cite{KLS}
\be
P_n(x; \beta, \gamma, \delta, N) = \kappa_n \:{_4}F_3 \left( {-n, -n+\beta-N, -s, s+\gamma+\delta+1 \atop -N, \beta+\delta+1, \gamma+1}; 1 \right) , \lab{Racah_P} \ee 
where
\be
x(s) = s(s+\gamma+\delta+1) \lab{x_s_Racah} \ee
The recurrence relation is
\be
P_{n+1}(x) +b_n P_n(x) +u_n P_{n-1}(x) = xP_n(x), \lab{rec_Racah} \ee
with
\be
b_n = -A_n -C_n, \quad u_n = A_{n-1}C_n, \lab{ub_Racah} \ee
where
\ba
&&A_n = \frac{(n+\beta-N)(n+\beta+\delta+1)(n+\gamma+1)(n-N)}{(2n+\beta-N)(2n+\beta-N+1)}, \nonumber \\
&& C_n= \frac{n(n+\beta)(n+\beta-\gamma-N-1)(n-\delta-N-1)}{(2n+\beta-N)(2n+\beta-N-1)}\lab{AC_Racah} \ea
The orthogonality of the Racah polynomials can be presented as \cite{KLS}
\be
\sum_{s=0}^N W_s P_n(x(s)) P_m(x(s)) = h_n \: \delta_{nm}, \lab{ort_Racah} \ee
where the weights are 
\be
W_s =M^{-1} \: \frac{(-N)_s(\beta+\delta+1)_s(\gamma+1)_s(\gamma+\delta+1)_s((\gamma+\delta+3)/2)_s}{s! (\gamma+\delta+2+N)_s(-\beta+\gamma+1)_s(\delta+1)_s ((\gamma+\delta+1)/2)_s},     \lab{W_racah} \ee
where the "total mass" is
\be
M= \frac{(-\beta)_N(\gamma+\delta+2)_N}{(1+\gamma-\beta)_N(\delta+1)_N}. \lab{total_M_Racah} \ee
The mirror-dual Rach polynomials $P_n^*(x)$ again belong to the Racah class and have the parameters
\be
\beta^* = -\beta, \: \gamma^* = \delta, \: \delta^* = \gamma . \lab{mirror_Racah} \ee

Necessary condition to get all the weights $W_s$ not depending on $s$ is to equate all parameters in the numerator Pochhammer symbols in \re{W_racah} to corresponding parameters in the denominator Pochhammer symbols. There are several solutions. However, up to the reflection $s \to N-s$ and shifting the linear grid, there are basically only two solutions:

\vspace{2mm}

(i) $\gamma=-1/2, \: \delta=1/2, \: \beta=N+1/2$. In this case $x_s = s(s+1)$.

\vspace{3mm}

(ii) $\gamma=\delta =-1/2, \: \beta= N+1/2$. In this case $x_s = s^2$.

\vspace{3mm}

The solution (ii) should be discarded because it leads to a singularity in the  the expression \re{W_racah}. This singularity may be avoided if one takes appropriate limiting procedure $\delta = -1/2 + \varepsilon, \: g = - 1/2 + \varepsilon$ with $\varepsilon \to 0$. This procedure yields the weights (up to a normalization factor)
\be
W_0 = 1, \: W_s  =2, \; s=1,2,\dots N \lab{W_0-nl} \ee  
which does not correspond to the Legendre case.

We thus arrived at the only possible type of the quadratic grid $x_s=s(s+1)$.

The recurrence coefficients are 
\be
u_n^* = \frac{n^2(2n-1)^2 \left((N+1)^2-n^2)\right)(2N+3-2n)(2N+1+2n)}{(4n+1)(4n-3)(4n-1)^2} \lab{u_Racah_2} \ee
and 
\be
b_n^* = \frac{(N+5/4)(N+3/4)}{8} \left( \frac{1}{4n-1} -  \frac{1}{4n-3}\right) + \frac{(N-n)(2n+2N+1)}{4} + \frac{3N}{4} + \frac{5}{32} \lab{b_Racah_2} \ee
Jacobi matrix $J$ corresponding to the Sturm sequence on the same quadratic grid $x_s=s(s+1)$ are obtained from the above coefficients by the mirror symmetry $b_n \to b_{N-n}, \: u_n \to u_{N+1-n}$. 

We thus have the
\begin{pr}
The Sturm sequence on quadratic grid $x_s=s(s+1)$  coincides with the sequence of the Racah polynomials $R_n(x;\alpha,\beta, \gamma, \delta)$ with the parameters   $\alpha= -N-1, \gamma=1/2, \: \delta=-1/2, \: \beta=-N-1/2$. 
\end{pr}

Consider now the quadratic grid with $\tau=2$, i.e. $x_s=s(s+2)$. If one puts 
\be
\delta= \gamma = 1/2, \quad \beta = N+3/2 \lab{par-rac-2} \ee
then the normalized Racah weights \re{W_racah} become 
\be
W_s = M^{-1} \: (x_s+1), \lab{W-Rac-2} \ee
where
\be
M=\frac{N(N+1)(2N+7)}{6}. \lab{M_racah-2} \ee

Let $P^{(0)}(x)$ be the orthogonal polynomials corresponding to the Legendre-type measure $W_s^{(0)}=(N+1)^{-1}$ on the grid $x_s=s(s+2)$. Then formula \re{W-Rac-2} indicates that the Racah  polynomials $P_n(x)$ with $\delta= \gamma = 1/2, \quad \beta = N+3/2$ are obtained from $P^{(0)}(x)$ by the Christoffel transform \cite{ZheR}
\be
P_n(x) = \frac{P_{n+1}^{(0)}(x) - V_n P_n^{(0)}(x)}{x+1}, \quad V_n = \frac{P_{n+1}^{(0)}(-1)}{P_n^{(0)}(-1)}. \lab{CT-Rac} \ee
Indeed, the Christoffel transform \cite{ZheR} is transformation of the orthogonal polynomials $P_n(x) \to \t P_n(x)$ which corresponds to the transformation of the orthogonality measure $d \mu(x)$
\be
d \tilde \mu(x)  = (x-a) d \mu(x) \lab{CT_measure} \ee
with some real parameter $a$ located outside the orthogonality interval.

Equivalently, the Christoffel transform can be presented as
\be
\t P_n(x) = \frac{P_{n+1}(x) - V_n P_n(x)}{x-a}, \quad V_n = \frac{P_{n+1}(a)}{P_n(a)}. \lab{CT-P} \ee 
The transformed polynomials $\t P_n(x)$ satisfy the recurrence relation 
\be
\t  P_{n+1}(x) + \t b_n \t P_n(x) + \t u_n \t P_{n-1}(x) = x \t P_n(x) \lab{CT_rec} \ee
with the transformed recurrence coefficients \cite{ZheR}
\be
\t u_n = u_n \: V_{n}/V_{n-1},  \quad \t b_n = b_{n+1} +V_{n+1}-V_n. \lab{CT_ub} \ee 
Conversely, the polynomials $P_n(x)$ can be obtained from the $q$-Hahn polynomials $Q_n(x)$ by the Geronimus transform \cite{ZheR}
\be
P_n(x) = \t P_n(x) - U_n \t P_{n-1}(x) , \lab{Ger} \ee
where
\be
U_n = \frac{\vphi_n}{\vphi_{n-1}}. \lab{G_phi} \ee
The sequence $\phi_n(z)$ is an arbitrary solution of the recurrence relation
\be
\vphi_{n+1} + \t b_n \vphi_n + \t u_n \vphi_{n-1} = a \vphi_n . \lab{vphi_def} \ee
We thus see that the polynomials of Legendre type on the quadratic grid $x_s=s(s+2)$ can be obtained from the Racah polynomials $\t P_n(x)$ on the same grid by the Geronimus transform  \re{Ger}, where the parameters of the Racah polynomials are $\alpha=-N-1, \: \delta= \gamma = 1/2, \: \quad \beta = N+3/2$.

In order to find corresponding Sturmian polynomials on the grid $x_s=s(s+2)$ we can use formulas   \re{ww_rel} and \re{mirror_Racah}. From these formula it is seen that the Sturmian polynomials are obtained from the Racah polynomials with the parameters $\alpha=_N-1, \: \beta=-N-3/2, \: \gamma=\delta=1/2$ by the Christoffel transformation \re{CT-P} with $a=-1$.

Hence we have 
\begin{pr}
The Sturm sequence $P_n(x), \: n=0,1,2\dots, N$ of the polynomials on the quadratic grid $x_s=s(s+2)$ coincides with polynomials obtained by the Christoffel transforms of the Racah polynomials $R_n(x) \equiv P_n(x; \alpha,\beta, \delta, \gamma)$ with parameters $\alpha=_N-1, \: \beta=-N-3/2, \: \gamma=\delta=1/2$:
\be
P_n(x) = \frac{R_{n+1}(x)- \frac{R_{n+1}(-1)}{R_n(-1)} R_n(x)}{x+1}. \lab{Chr_Rac} \ee
\end{pr}
These polynomials do not belong to the Askey scheme \cite{KLS} and hence they do not satisfy any second-order difference equation on the grid $x_s=s(s+2)$.

\section{Exponential grid}
\setcounter{equation}{0}
The exponential grid can be reduced by affine transformations to the canonical form
\be
x_s = q^{-s}, \quad s=0,1,\dots N . \lab{x_exponent} \ee 
We assume that $0<q<1$.

In order to find a solution,  it is natural to consider the q-Hahn polynomials.

These polynomials have the explicit expression \cite{KLS}
\be
Q_n(x; \alpha,\beta,N) = {_3}\phi_2 \left( {q^{-n}, \alpha \beta q^{n+1} , x \atop q^{-N}, \alpha q }  \mid q;q  \right) . \lab{q-Hahn_hypo} \ee
They are orthogonal on the exponential grid \re{x_exponent}
\be
\sum_{s=0}^N W_s Q_n(q^{-s}) Q_m(q^{-s}) = h_n \: \delta_{nm} \lab{ort-q-hahn} \ee
with 
\be
W_s= M^{-1} \:\frac{(\alpha q;q)_s(q^{-N};q)_s}{(q;q)_s (\beta^{-1} q^{-N};q)_s} (\alpha\beta q)^{-s} \lab{W_q-hahn} \ee
and
\be
M= \frac{(\alpha \beta q^2;q)_N}{(\beta q;q)_N(\alpha q)^N}, \lab{mass_q-hahn} \ee
where $(x;q)_s$ stands for the q-Pochhammer symbol
\[
(x;q)_s =(1-x) (1-xq) \dots (1-xq^{s-1}).
\]
The q-Hahn polynomials satisfy the three-term recurrence relation
\be
Q_{n+1}(x) + b_n Q_n Q_n(x) + u_n Q_{n-1}(x) = xQ_n(x) \lab{rec_Q} \ee 
with \cite{KLS}
\be
u_n = A_{n-1} C_n, \quad b_n = 1-A_n-C_n, \lab{recQ_AC} \ee
where
\be
A_n ={\frac { \left( 1-{q}^{n-N} \right)  \left( 1-a{q}^{n+1} \right) 
 \left( 1-ab{q}^{n+1} \right) }{ \left( 1-ab{q}^{2\,n+1} \right) 
 \left( 1-ab{q}^{2\,n+2} \right) }} \lab{A_q-hahn} \ee
and 
\be
C_n = -{\frac {a{q}^{n-N} \left( 1-{q}^{n} \right)  \left( 1-b{q}^{n}
 \right)  \left( 1-ab{q}^{n+N+1} \right) }{ \left( 1-ab{q}^{2\,n+1}
 \right)  \left( 1-ab{q}^{2\,n} \right) }}. \lab{C_q-hahn} \ee

In contrast to the case of the ordinary Hahn polynomials, it is impossible to obtain the Legendre-type weights $W_s =(N+1)^{-1}$ for q-Hahn polynomials .

Nevertheless, we can achieve the desired result by putting $\alpha=\beta=1$. Then the weights become
\be
W_s = M^{-1} q^{-s}. \lab{W_exp} \ee
This means that the measure \re{W_exp} is obtained from the uniform Legendre-type measure by multiplying to the argument $x$ which takes the values $x(s)=q^{-s}$ on the grid $x(s)$. 

Let $P^{(0)}(x)$ be the orthogonal polynomials corresponding to the Legendre-type measure $W_s^{(0)}=(N+1)^{-1}$ on the grid $q^{-s}$. Then formula \re{W_exp} indicates that the q-Hahn polynomials $Q_n(x)$ with $\alpha=\beta=1$ are obtained from $P^{(0)}(x)$ by the Christoffel transform \cite{ZheR}
\be
Q_n(x) = \frac{P_{n+1}^{(0)}(x) - V_n P_n^{(0)}(x)}{x}, \quad V_n = \frac{P_{n+1}^{(0)}(0)}{P_n^{(0)}(0)}. \lab{CT-q} \ee
Conversely, the polynomials $P_n^{(0)}(x)$ can be obtained from the $q$-Hahn polynomials $Q_n(x)$ by the Uvarov transform \cite{Uvarov}, which can be considered as a special case of the Geronimus transform \cite{ZheR}
\be
P_n^{(0)}(x) = Q_n(x) - U_n Q_{n-1}(x) , \lab{Uvarov} \ee
where
\be
U_n = \frac{F_n(0)}{F_{n-1}(0)}. \lab{U_F} \ee
The functions $F_n(z)$ are defined as 
\be
F_n(z)  = \int \frac{Q_n(x)}{z-x} d \mu(x) .  \lab{F_n_def} \ee
These functions are known to yield the second linear independent solution of the same recurrence relation as for the orthogonal polynomials $P_n(x)$ \cite{NSU}
\be
F_{n+1}(z) + b_n F_n(z) + u_n F_{n-1}(z) = z F_n(z). \lab{rec_F} \ee
In our concrete case these functions can be calculated as
\be
F_n(0)= \sum_{s=0}^N \frac{W_s Q_n(q^{-s})}{-q^{-s}} = \kappa \: \sum_{s=0}^N Q_n(q^{-s}) \lab{F_n(0)} \ee
with some normalization constant $\kappa$ which does not depend on $n$.

We thus have the
\begin{pr}
The Legendre type polynomials $P_n^{(0)}(x)$ on the exponential grid $x_s=q^{-s}$ are given by the Uvarov transform \re{Uvarov}-\re{U_F} of the q-Hahn polynomials $Q_n(x;\alpha,\beta,N)$ with $\alpha=\beta=1$, where $F_n(0)$ has expression \re{F_n(0)}. 
\end{pr}
In order to obtain explicit expression for the Sturm polynomials on the exponential grid we notice that the reflection  $Q_n(x) \to Q_n^*(x)$ for q-Hahn polynomials is equivalent to the transformation
\be
a^* = b^{-1} q^{-N-1}, \quad b^* = a^{-1} q^{-N-1} . \lab{*-tr-q-hahn} \ee
Moreover, from the formula \re{ww_rel} we see that the Uvarov transform of the polynomials $P_n^*(x)$ corresponds to the Christoffel transform of the polynomials $P_n(x)$. Omitting obvious technical details, we arrive at 
\begin{pr}
The Sturm sequence $P_n(x), \: n=0,1,2\dots, N$ of the polynomials on the exponential grid $x_s=q^{-s}$ coincides with polynomials obtained by the Christoffel transforms of the q-Hahn polynomials with parameters $a=b=q^{-N-1}$:
\be
P_n(x) = \frac{Q_{n+1}(x)- \frac{Q_{n+1}(0)}{Q_n(0)} Q_n(x)}{x}, \lab{Chr_Q} \ee
where $Q_n(x) \equiv Q_n(x; q^{-N-1}, q^{-N-1},N)$.
\end{pr} 
As a direct consequence of this proposition, using \re{CT_ub}, we derive explicit   expressions for the recurrence coefficients of the Sturm polynomials:
\be
u_n = u_n^{(0)} V_{n}/V_{n-1}, \quad b_n = b_{n+1}^{(0)} +V_{n+1}-V_n, \lab{ub-q-Sturm} \ee 
where 
\be
V_n = \frac{Q_{n+1}(0)}{Q_n(0)} \lab{V_n_0} \ee
and where $u_n^{(0)}$ and $b_n^{(0)}$ are recurrence coefficients \re{recQ_AC} of the q-Hahn polynomials with the parameters $a=b=q^{-N-1}$.

Note that in contrast to the case of linear grid, the polynomials \re{Chr_Q} do not belong to the Askey scheme. Hence they don't satisfy any second order difference eigenvalue equation.

\section{Trigonometric grids}
\setcounter{equation}{0}
In this section we consider two special cases of the trigonometric grids connected with roots of unity.

The first grid is given by
\be
x_s= -\cos\left(\omega (s+1/2)\right), \quad s=0,1,2,\dots, N, \lab{roots_T} \ee
where
\be
\omega = \frac{\pi}{N+1}. \lab{omega_trig} \ee
This trigonometric grid is related to the Chebyshev polynomials  polynomials of the first kind $T_n(x)$ \cite{MH}. Indeed, these (monic) polynomials are defined as
\be
T_0=1, \quad  T_n(x) = 2^{1-n} \: \cos(\theta n), \quad n=1,2,3,\dots  \lab{P_T} \ee
with 
\be
x=\cos \theta \lab{x_cos} \ee
They satisfy the recurrence relation
\be
T_{n+1}(x) + u_n T_{n-1}(x) = xT_n(x) , \quad n=1,2,3, \dots \lab{rec_T} \ee
where the recurrence coefficients are
\be
u_1 = 1/2, \; u_n = 1/4, \; n=2,3,\dots \lab{u_T} \ee
The monic Chebyshev polynomials of the second kind are defined as
\be
U_n(x) = 2^{-n} \: \frac{\sin (\theta(n+1))}{\sin \theta}, \quad n=0,1,2,\dots \lab{P_U} \ee
They satisfy the same recurrence relation
\be
U_{n+1}(x) + u_n U_{n-1}(x) = xU_n(x) , \quad n=1,2,3, \dots \lab{rec_U} \ee
with $u_n=1/4$ for all $n=1,2,\dots$. Thus the only difference between recurrence relations for polynomials $T_n(x)$ and $U_n(x)$ is expression of the first recurrence coefficient: $u_1=1/2$ for $T_n(x)$ and $u_1=1/4$ for $U_n(x)$.

There is obvious relation between these polynomials \cite{MH}:
\be
T_{n+1}'(x) = (n+1) U_n(x), \quad n=0,1,2,\dots  \lab{T_U} \ee
The roots $x_s$ of the Chebyshev polynomial $T_{N+1}(x)$ coincide with \re{roots_T}. On these roots the Chebyshev polynomials $T_n(x)$ satisfy the finite orthogonality relation \cite{MH}
\be
\sum_{s=0}^N T_n(x_s) T_m(x_s) =(N+1) h_n \delta_{nm}. \lab{ort_T} \ee
Hence  the Chebyshev polynomials of the first kind $T_n(x)$ satisfy orthogonality relation  of Legendre type. This allows one to construct corresponding Sturmian orthogonal polynomials.

Indeed, let us start with prescribed grid \re{roots_T}. Then the initial Sturmian polynomial $P_{N+1}(x)$ coincides with the Chebyshev polynomial having zeros on this grid:  $P_{N+1}(x)=T_{N+1}(x)$. Using formula \re{T_U}, we conclude that the companion Sturmian polynomial $P_N(x)$ coincides with the Chebushev polynomial of the second kind: $P_N(x)=(N+1)^{-1} P_{N+1}'(x)= U_{N}(x)$. Applying the Sturm algorithm, we then reconstruct the whole chain of Sturmian orthogonal polynomials $P_{N-1}(x), P_{N-2}(x), \dots, P_1(x)$. 

In order to recognize these polynomials $P_n(x)$, we can use Proposition {\bf 1}. The polynomials $P_n(x)$ are mirror-dual with respect to the polynomials $T_0,T_1(x), T_2(x), \dots T_{N}(x)$. Hence we have the recurrence coefficients for them
\be
b_n =0, \quad n=0,1,2, \dots, N, \quad u_1 =u_2= \dots = u_{N-1}= 1/4, \: u_{N}=1/2 . \lab{rec_mT} \ee   
Comparing recurrence coefficients \re{rec_mT} we arrive at 
\begin{pr}
For the trigonometric grid \re{roots_T} the two initial Sturmian polynomials are $P_{N+1}(x) = T_{N+1}(x)$ and $P_N(x)=U_N(x)$. All further orthogonal polynomials coincide with the Chebyshev polynomials of the second kind: $P_n(x) = U_n(x), \; n=0,1,\dots, N-1$. These polynomials satisfy the recurrence relation 
\be
P_{n+1}(x) + u_n P_{n-1}(x) = xP_n(x), \quad n=1,2,\dots, N \lab{rec_trig} \ee
with the recurrence coefficients $u_n$ given by \re{rec_mT}.
\end{pr}
It is easy to find the orthogonality relation for the above polynomials $P_n(x)$
\be
\sum_{s=0}^N P_n(x_s) P_m(x_s) w_s  = h_n \: \delta_{nm}. \lab{ort_S_trig} \ee
Indeed, we already know that the mirror-dual orthogonality relation is of Legendre type \re{ort_T}. Hence, by \re{ww_rel} we have
\be
w_s = \frac{2}{N+1} \: \left(\frac{\sin \theta_s }{\sin (\theta_s (N+1))}\right)^2 =\frac{2}{N+1} \: \sin^2 \theta_s , \lab{w_trig_norm} \ee
where
\be
\theta_s =  \frac{\pi (s+1/2)}{N+1}. \lab{theta_s} \ee
These weights are normalized, i.e.
\be
\sum_{s=0}^N w_s =1 . \lab{norm_trig} \ee
As expected, the weights \re{w_trig_norm} are well known discrete orthogonality weights for the Chebyshev polynomials $U_n(x)$ on the roots of the Chebyshev polynomial $T_{N+1}(x)$ \cite{MH}.

The second case of the trigonometric grid coincides with zeros of the Chebyshev polynomial $U_{N+1}(x)$ of the second kind:
\be
x_s = -\cos\left(\frac{\pi(s+1)}{N+2}\right), \quad s=0,1, \dots, N . \lab{U_grid} \ee
In order to describe the corresponding Sturmian sequence of polynomials $P_n(x), \: x=0,1,\dots, N$ we recall the well known properties of the ultrasperical polynomials. The monic ultraspherical polynomials  $C_n^{(\lambda)}(x)$ depend on one parameter $\lambda$ and are defined as \cite{KLS}
\be
C^{(\lambda)}(x) = \kappa_n \: {_2}F_1 \left( {-n , n+2 \lambda \atop \lambda+ 1/2 } \left| \frac{1-x}{2} \right . \right),  \lab{C_def} \ee
where the normalization coefficient is 
\be
\kappa_n = 2^{-n} \: \frac{(2 \lambda)_n}{(\lambda)_n}. \lab{nor_C} \ee
They satisfy the 3-term recurrence relation \cite{KLS}
\be
C_{n+1}^{(\lambda)}(x) + u_n C_{n-1}^{(\lambda)}(x) = x C_{n}^{(\lambda)}(x) \lab{rec_C} \ee  
with 
\be
u_n = \frac{n(n+2 \lambda-1)}{4(n+\lambda)(n+\lambda-1)}. \lab{rec_u_C} \ee
The Chebyshev polynomials are special cases of the ultrapsherical polynomials
\be
T_n(x) = C_n^{(0)}(x), \quad U_n(x) = C_n^{(1)}(x). \lab{TU_C} \ee
Moreover, there is simple formula \cite{KLS}
\be
\frac{d C_n^{(\lambda)}(x)}{d x} = n C_{n-1}^{(\lambda+1)}(x) \lab{diff_C} \ee
from which it follows that
\be
\frac{d U_{n+1}(x)}{d x} = (n+1) \: C_n^{(2)}(x) . \lab{diff_U_C} \ee
We can now take derivatives from all members of the Sturmian sequence $U_0, U_1(x), \dots, U_N(x), T_{N+1}(x)$ corresponding to the previous type of trigonometric grid. We then have the sequence $$C_0^{(2)}(x),C_1^{(2)}(x), C_2^{(2)}(x), \dots, C_N^{(2)}(x), U_{N+1}(x),$$ 
where all polynomials apart from the final one are the ultrapsherical polynomials $C_n^{(2)}(x)$ and the final polynomial $P_{N+1}(x)$ coincides with the Chebyshev polynomial $U_{N+1}(x)$. The roots of the polynomial $P_{N+1}(x)$ are described by \re{U_grid}. Moreover, due to \re{diff_U_C} we have $(N+1) P_N(x) = P_{N+1}'(x)$. Hence we have the needed Sturm sequence of the polynomials $P_n(x)$ satisfying the recurrence relation 
\be
P_{n+1}(x) + u_n P_{n-1}(x) = x P_{n}(x) \lab{rec_St2} \ee
with the recurrence coefficients
\be
u_n = \frac{n(n+3)}{4(n+2)(n+1)} , \quad  n=1,2,\dots, N-1, \quad u_N = \frac{N}{2(N+1)}. \lab{rec_u_St2} \ee
We thus see that both trigonometric grids \re{roots_T} and \re{U_grid} correspond to elementary solutions of the Sturm sequence related with the Chebyshev polynomials $T_n(x)$ and $U_n(x)$.

\vspace{5mm}

\section{Conclusion}
\setcounter{equation}{0}
We have derived explicit systems of finite orthogonal polynomials corresponding to classical Sturm problem for the pair of polynomials $P(x)$ and $P'(x)$ on prescribed grids: linear, exponential and for special types of quadratic and trigonometric grids. In all the above cases, the corresponding polynomials are related to classical polynomials from the Askey scheme \cite{KLS}.

The open problems are constructing of Sturmian polynomials on other types of grids: Askey-Wilson and Bannai-Ito. Another interesting problem is related to the quadratic and trigonometric grid. We were able to construct explicitly the Sturm polynomials for the grid $x_s=s(s+1)$ and $x_s=s(s+2)$. What about generic quadratic grid $x_s=s(s+\tau)$ with an arbitrary parameter $\tau$?  And what are other examples of trigonometric grids admitting explicit solution for the Sturm sequence?

Another open problem arises in connection with Sturmian polynomials on the quadratic $x_s=s(s+2)$ and exponential $x_s=q^{-s}$ grids. As we already mentioned, these polynomials are obtained from the classical ones (i.e. from the Racah and q-Hahn polynomials) by the Christoffel transform. This means that they are not "classical" in sense of classification of the Aseky scheme \cite{KLS}. In particular, they don't satisfy a linear difference equation of second order. On the other hand, it is known that some Christoffel transformed classical polynomials of discrete variable do satisfy difference equations of {\it higher} orders \cite{Duran}. One can ask whether or not the Sturmian polynomials on the above lattices satisfy some difference equations of order 4 or more.

\bigskip\bigskip
{\Large\bf Acknowledgments}
The work is supported by the National
Science Foundation of China (Grant No.11771015). The author thanks to Luc Vinet for discussion.

\vspace{5mm}

\bb{99}

\bi{Beardon} A.F.Beardon, {\it Sums of Powers of Integers}, Am.Math.Month. 1996, 201--213.

\bi{BS} C. de Boor, E.B. Saff, {\it Finite sequences of orthogonal polynomials connected by a Jacobi matrix}, Linear Algebra
Appl. {\bf 75} (1986) 43--55.

\bi{Borodin}  A. Borodin, {\it Duality of orthogonal polynomials on a finite set}, J. Statist. Phys. {\bf 109} (2002) 1109--1120.

\bi{Chi} T. Chihara, 
{\it An Introduction to Orthogonal Polynomials}, 
Gordon and Breach, NY, 1978.

\bi{DS} H. Dette, W.J. Studden, {\it On a new characterization of the classical orthogonal polynomials}, J. Approx. Theory {\bf 71}
(1) (1992) 3--17.

\bi{Duran} A. J. Dur\'an, {\it Orthogonal polynomials satisfying higher-order difference
equations}, Constr. Approx. {\bf 36} (2012), 459--486.

\bi{persym} V.Genest, S. Tsujimoto, L.Vinet and A.Zhedanov, {\it Persymmetric Jacobi matrices, isospectral deformations and orthogonal polynomials}. J. Math. Anal. Appl. {\bf 450} (2017), no. 2, 915--928.

\bibitem{KLS} R. Koekoek, P.A. Lesky, and R.F. Swarttouw. {\it Hypergeometric orthogonal polynomials and their q-analogues}. Springer, 1-st edition, 2010.

\bi{MH} J.C.Mason, D.C.Handscomb, {\it Chebyshev Polynomials}, Chapman \& Hall/CRC, 2003.

\bi{NSU} A.F. Nikiforov, S.K. Suslov, and V.B. Uvarov, 
{\it Classical Orthogonal Polynomials of a Discrete Variable},
Springer, Berlin, 1991.

\bi{Prasolov} V.Prasolov, {\it Polynomials}, (Translated from Russian), Springer-Verlag Berlin Heidelberg, 2004.

\bi{SZ_Duke} V.Spiridonov and A.Zhedanov, {\it Zeros and Orthogonality of the
Askey-Wilson Polynomials for  q a root of unity}, {\bf 89} (1997),  283--305.

\bi{Syl} J.J.Sylvester, {\it On a theory of the syzygetic relations of two rational integral functions, comprising an application to the theory of Sturm's functions, and that of the greatest algebraical common measure}. Phil. Trans. R. Soc. Lond. {\bf 143},  407--548.

\bi{Ter} P.Terwilliger, {\it Two linear transformations each tridiagonal with respect to an eigenbasis of
the other}, Lin.Alg.Appl. {\bf 330} (2001), 149--203.

\bi{Uvarov}  V.B. Uvarov, {\it The connection between systems of polynomials that are orthogonal with respect to different distribution
functions}, Zh.Vychisl.Matem. i Mat.Fiz. {\bf 9} (1969) 1253--1262 (in Russian).

\bi{VZ_Bochner} L.Vinet and A.Zhedanov, {\it Generalized Bochner theorem: Characterization of the Askey-Wilson polynomials}, J.Comp.Appl.Math. {\bf 211} (2008) 45--56.

\bi{VZ_dual} L.Vinet and A.Zhedanov, {\it A characterization of classical and semiclassical orthogonal polynomials from their dual polynomials}, J.Comp.Appl.Math.,  {\bf 172} (2004) 41--48.

\bi{Waerden} Van der Waerden, Algebra, vol. I, Springer, 1991.

\bi{ZheR} A.Zhedanov, {\it Rational spectral transformations and orthogonal polynomials}, J.Comp.Appl.Math. {\bf 85} (1997) 67--86.

\end{thebibliography}

\end{document}